\documentclass[12pt]{amsart}
\usepackage{amssymb,amsmath}
\def\contrazione{\raisebox{1pt}{\,{\mbox{\tiny{$|\!\raisebox{-0.7pt}{\underline{\hphantom{X}}}$}}}\,}}

\def\Om{\Omega}

\def\NO#1{||#1||^2}
\def\no#1{||#1||}

\numberwithin{equation}{section}
\def\db{\bar\partial}
\def\db*{\bar\partial^*}
\def\T{\text}

\def\simleq{\underset\sim<}
\def\1#1{\overline{#1}}
\def\2#1{\widetilde{#1}}
\def\3#1{\widehat{#1}}
\def\4#1{\mathbb{#1}}

\def\5#1{\frak{#1}}
\def\6#1{{\mathcal{#1}}}
\def\C{{\4C}}
\def\R{{\4R}}

\def\Z{{\4Z}}
\def\sumK{\underset{|K|=k-1}{{\sum}'}} 
\def\sumJ{\underset{|J|=k}{{\sum}'}}
\def\sumij{\underset {ij=1,...,n-1}{{\sum}}}

\def\sumj{\underset {j=1,\dots,n-1}{{\sum}}}

\def\di{\partial}
\def\dib{\bar\partial}
\begin{document}
\abstract
We prove local hypoellipticity of the complex Laplacian $\Box$ and of the Kohn Laplacian $\Box_b$  in a pseudoconvex  boundary when,  for a system of cut-off $\eta$, the gradient $\di_b\eta$ and the Levi form $\di_b\dib_b\eta^2$ are subelliptic multipliers in the sense of \cite{K79}. 
\newline
MSC: 32F10, 32F20, 32N15, 32T25 
\endabstract
\title[Hypoellipticity of the $\bar\partial$-Neumann problem...]{Hypoellipticity of the $\bar\partial$-Neumann problem by means of subelliptic multipliers}
\author[L.~Baracco, S.~Pinton and G.~Zampieri]{Luca Baracco, Stefano Pinton and  Giuseppe Zampieri}
\address{Dipartimento di Matematica, Universit\`a di Padova, via 
Trieste 63, 35121 Padova, Italy}
\email{baracco@math.unipd.it,pinton@math.unipd.it,  
zampieri@math.unipd.it}
%\subjclass{}
\maketitle
%\tableofcontents
% Standard sets
\def\Giialpha{\mathcal G^{i,i\alpha}}
\def\cn{{\C^n}}
\def\cnn{{\C^{n'}}}
\def\ocn{\2{\C^n}}
\def\ocnn{\2{\C^{n'}}}
% Abbreviations
\def\const{{\rm const}}
\def\rk{{\rm rank\,}}
\def\id{{\sf id}}
\def\aut{{\sf aut}}
\def\Aut{{\sf Aut}}
\def\CR{{\rm CR}}
\def\GL{{\sf GL}}
\def\Re{{\sf Re}\,}
\def\Im{{\sf Im}\,}
\def\codim{{\rm codim}}
\def\crd{\dim_{{\rm CR}}}
\def\crc{{\rm codim_{CR}}}
\def\phi{\varphi}
\def\eps{\varepsilon}
\def\d{\partial}
\def\a{\alpha}
\def\b{\beta}
\def\g{\gamma}
\def\G{\Gamma}
\def\D{\Delta}
\def\Om{\Omega}
\def\k{\kappa}
\def\l{\lambda}
\def\L{\Lambda}
\def\z{{\bar z}}
\def\w{{\bar w}}
\def\Z{{\1Z}}
\def\t{{\tau}}
\def\th{\theta}
\emergencystretch15pt
\frenchspacing
\newtheorem{Thm}{Theorem}[section]
\newtheorem{Cor}[Thm]{Corollary}
\newtheorem{Pro}[Thm]{Proposition}
\newtheorem{Lem}[Thm]{Lemma}
\theoremstyle{definition}\newtheorem{Def}[Thm]{Definition}
\theoremstyle{remark}
\newtheorem{Rem}[Thm]{Remark}
\newtheorem{Exa}[Thm]{Example}
\newtheorem{Exs}[Thm]{Examples}
\def\Label#1{\label{#1}}
\def\bl{\begin{Lem}}
\def\el{\end{Lem}}
\def\bp{\begin{Pro}}
\def\ep{\end{Pro}}
\def\bt{\begin{Thm}}
\def\et{\end{Thm}}
\def\bc{\begin{Cor}}
\def\ec{\end{Cor}}
\def\bd{\begin{Def}}
\def\ed{\end{Def}}
\def\br{\begin{Rem}}
\def\er{\end{Rem}}
\def\be{\begin{Exa}}
\def\ee{\end{Exa}}
\def\bpf{\begin{proof}}
\def\epf{\end{proof}}
\def\ben{\begin{enumerate}}
\def\een{\end{enumerate}}
\def\dotgamma{\Gamma}
\def\dothatgamma{ {\hat\Gamma}}

\def\simto{\overset\sim\to\to}
\def\1alpha{[\frac1\alpha]}
\def\T{\text}
\def\R{{\Bbb R}}
\def\I{{\Bbb I}}
\def\C{{\Bbb C}}
\def\Z{{\Bbb Z}}
\def\Fialpha{{\mathcal F^{i,\alpha}}}
\def\Fiialpha{{\mathcal F^{i,i\alpha}}}
\def\Figamma{{\mathcal F^{i,\gamma}}}
\def\Real{\Re}
%
%
%
%\endtopmatter
\section{Introduction}
For a pseudoconvex domain $\Om\subset\subset\C^n$ with $C^\infty$-boundary $b\Om$, we consider the problem of the local regularity of the canonical solution of $\dib_b$ and of the $\dib$-Neumann problem at a point $z_o\in b\Om$. We form the Kohn Laplacian $\Box_b=\dib^*_b\dib_b+\dib_b\dib^*_b$ and the complex Laplacian $\Box= \dib^*\dib+\dib\dib^*$. The first problem can be restated in terms of the hypoellipticity of $\Box_b$
$$
\T{\it( Hypoellipticity)}\quad\Box_bu\in C^\infty_{z_o}\quad\T{implies}\quad u\in C^\infty_{z_o}.
$$
In the same way the hypoellipticity of $\Box$ is defined. We search for general criteria of hypoellipticity. It was firstly noticed by Kohn that the presence in supporting complex hypersurfaces of propagators   of boundary smoothness  of holomorphic functions  prevents from hypoellipticity.
A related phenomenon is that of the propagation of holomorphic extendibility. According to \cite{HT83}, this takes place along complex curves. However, in the exponentially degenerate case, it was proved by \cite{BKZ12} that real curves are also propagators. This is the case of the lines $\R_{y_j}$ for the tube domain $2x_2=e^{-\frac1{\sum_{j=1}^{n-1}|x_j|^s}}$ for $s\ge1$. This propagation matches the non-hypoellipticity of $\Box_b$ proved by \cite{Ch02} in $\C^2$. Instead, if $s<1$, the argument for propagation of \cite{BKZ12} breacks down; again, this is in accordance with the hypoellipticity which occurs  as a consequence of ``superlogarithmic" estimates (cf. this section below). Thus propagation and hypoellipticity appear opposite one to another. 

As for classical positive results on hypoellipticity, we recall that this is generally obtained through estimates on forms $v$ of degree $k\in[1,n-2]$ such as 
\begin{equation}
\Label{subelliptic}
\T{\it (Subelliptic)}\quad\no{v}_\epsilon\simleq \no{\dib_bv}+\no{\dib^*_bv},
\end{equation}
or
\begin{equation}
\Label{suplog}
\T{\it (Superlogarithmic)}\quad \no{\log(\Lambda)v}\simleq \delta(\no{\dib_bv}+\no{\dib^*_bv})+c_\delta\no{v}_{-1},
\end{equation}
for any $\delta$ and for suitable $c_\delta$. Models are ``decoupled" domains $2x_n=\sum_{j=1}^{n-1}h^j(z_j)$ with
$$
\T{\it (Subelliptic)}\quad h^j=|z_j|^{2m_j}\quad\T{or}\quad h^j=x_j^{2m_j},
$$
$$
\T{\it (Superlogarithmic)}\quad h^j=e^{-\frac1{|z_j|^s}}\quad\T{or}\quad h^j=e^{-\frac1{|x_j|^s}},\quad s<1,
$$
where we can replace the power $|z_j|^s$, $s<1$ by $|z_j|\log|z_j|$ and similarly for $x_j$ (cf. \cite{KZ10}).
To get hypoellipticity from  \eqref{subelliptic} (cf. \cite{KN65}), one substitutes $\eta\Lambda^su$ for $v$ where $\eta$ ranges in a system of cut-off and $\Lambda^s$ is the standard tangential elliptic operator of order $s$. The problem is to control the commutators $[\dib^{(*)}_b,\eta\Lambda^s]$. First, these are estimated by $|\di\eta|\Lambda^s+c_s\Lambda^s$; next, one controls $c_s$ by a small constant produced by Sobolev interpolation,
 $|\di\eta|\Lambda^s$ by induction, and gets
\begin{equation}
\Label{apriori}
\no{\eta u}_s\simleq \no{\eta'\dib_bu}_s+\no{\eta'\dib^*_bu}_s+\no{u}_0\quad\T{for $\eta'\succ\eta$ i.e. $\eta'|_{\T{supp}\,\eta}\equiv1$},
\end{equation}
which is sufficient for hypoellipticity. To get the same conclusion \eqref{apriori} starting  from \eqref{suplog}, one replaces $\Lambda^s$ by the pseudodifferential operator $R^s$ with symbol $\sigma(R^s)=\Lambda_\xi^{s\sigma(z)}$ for $\eta\prec\sigma\prec\eta'$ and notices that $\eta\Lambda^s\prec \eta'R^s+O(\Lambda^{-\infty})$, $|\di\eta| R^s=O(\Lambda^{-\infty})$, 
$|[\dib_b^{(*)},R^s]|\le c_s\log(\Lambda)R^s$ and controls $c_s<<\delta^{-1}$ where $\delta$ is the small constant in \eqref{suplog} (cf. \cite{K02}).

But hypoellipticity is not entirely ruled by estimates. In \cite{K00}, Kohn proves hypoellipticity for boundaries defined by $2x_n=h(z',y_n)$ such that
\begin{itemize}
\item[(i)] there are subelliptic estimates for $|z'|\neq0$,
\item[(ii)] $h_{\bar z_j}$ are subelliptic multipliers.
\end{itemize}
In this situation, taking a cut-off $\chi$ of one real variable and setting $\zeta=\Pi_{j=1,...,n-1}\chi(|z_j|)$, $\theta=\chi(|y_n|)$, $\eta=\zeta\theta$, and denoting by $\bar L_j,\,\,j=1,...,n-1$ a system of $(0,1)$ vector fields, we have
$$
[\bar L_j,\eta]=\underset{\T{controlled by (i)}}{\underbrace{\zeta_{z_j}\theta}}+\underset{\T{controlled by (ii)}}{\underbrace{h^j_{\bar z_j}\zeta\dot\theta}}.
$$
The model  is 
$$
2x_n=e^{-\frac1{(\sum_j|z_j|)^s}}\quad\T{for any $s>0$}.
$$
When $s<1$, hypoellipticity was already obtained from superlogarithmicity even with $z_j$ replaced by $x_j$; when $s\ge1$, the conclusion is new and does not hold for $x_j$ (cf. \cite{BKZ12} and \cite{Ch02} already mentioned above). It remained open the problem of the hypoellipticity of domains with model
\begin{equation}
\Label{pinton}
2x_n=\sum_je^{-\frac1{|z_j|^s}},\quad s\ge1,
\end{equation}
in which summation is not taken at exponent. In this case (i) and (ii) do not hold at the points of the ``cross" $z_j=0$ for some $j=1,...,n-1$. A first answer to this question has been given in \cite{BPZ13} where hypoellipticity is stated on a class of domains which contains \eqref{pinton}. This is obtained by modifying the localized  ``bad" vector field $\zeta\theta T$ into
$$
(T)_{\zeta\theta}:=\zeta\theta T-\Big(\sum_j (h^j_{z_j\bar z_j})^{-1}L_j(\zeta\theta)\bar L_j+(h^j_{z_j\bar z_j})^{-1}\bar L_j(\zeta\theta) L_j\Big)
$$
(cf. \cite{DT95}). The class of domains in question is that for which the coefficients $(h^j_{z_j\bar z_j})^{-1}\overset{(-)}L_j(\zeta\theta)$ are well defined, that is, the zeroe's of $\overset{(-)}L_j(\zeta\theta)$ balance those of $h^j_{z_j\bar z_j}$. 

In the present paper, we give the geometric solution to the problem. Hypoellipticity holds whenever 
\begin{equation}
\Label{nomen}
\T{$\dib_b\eta$ and $\di_b\dib_b\eta^2$ are subelliptic multipliers,}
\end{equation}
over ``positively microlocalized" forms $u^+$. 
The model  is 
$$
2x_n=\sum_j e^{-\frac1{|z_j|^{s_j}}}x_j^{2m_j}\quad\T{any $s_j>0$ and $m_j\ge0$.}
$$
The idea of the proof is to insert the cut-off $\eta$ into the weight $e^{-\phi}$,  $\phi=-\log\eta^2+t|z|^2,\,\,z'\in T^\C b\Om$, which occurs in the ``basic estimate".  
This dispenses from controlling $[\dib_b^{(*)},\eta]$ and reduces the problem only  to the error in the Levi form and in the adjunction (in addition  to the commutator $[\dib^{(*)}_b,\Lambda^s]$, as usual):
\begin{equation*}
\begin{cases}
e^{-\phi}\Big(\di_b\dib_b\phi-\di_b\dib_b(t|z|^2)\Big)\sim\di_b\dib_b\eta^2,
\\
e^{-\phi}\Big((\dib^*_\phi)_b-(\dib^*_{t|z'|^2})_b\Big)\sim\di_b\eta.
\end{cases}
\end{equation*}
Thus, by the aid of \eqref{nomen}, the basic estimate turns into a regularity estimate with cut-off. Note that
 the single entries of $\di_b\dib_b\eta^2$ and $\di\eta$ need not to be subelliptic multipliers for all components of $u$ but just for those that they ``pick up".
\vskip0.3cm
\noindent
{\it Acknowledgments. } The paper was accomplished at Sao Paulo USP in November 2013. The authors are grateful to Paulo Domingo Cordaro for friendly hospitality and fruitful discussions.

%The proof is obtained from ``twisted" basic estimates with norm and adjunction weighted by $\eta$ (in addition to $e^{-t|z'|^2}$). In this way, one can see (cf. %Proposition~\ref{p1.1} below) that what really matters are not the commuatators $[\bar L_j, \eta]u$ themselves, but more precisely $\dib_b\eta\contrazione u$ and $\di_b\dib_b\eta(u,u)$.

\section{The main result}
Let $\Om\subset\subset\C^n$ be a smooth pseudoconvex domain and $z_o=0$ a boundary point.
\bt
\Label{t1.1}
Assume that there is a system of smooth cut-off $\eta$ in a neighborhood of $0$ such that
\begin{multline}
\Label{1.1}
\T{ $\di_b\eta$ and $\di_b\dib_b\eta^2$ are subelliptic multipliers in positive microlocalization}
\\
\T{ in any degree $k\in[1,n-1]$ (cf. \cite{K79})}.
\end{multline}
Then $\Box_b$ and $\Box$ are $C^\infty$-hypoelliptic at $0$.
\et
The main tool in the proof is the proposition below. Let $\mathcal H_b=\ker \Box_b$ be the space of harmonic forms.
\bp
\Label{p1.1}
Assume that for a system of cut-off $\eta$, \eqref{1.1} is satisfied. Then for any $\eta$  and for suitable $\eta'\succ\eta$, that is $\eta'|_{\T{supp}\,\eta}\equiv1$, we have
\begin{multline}
\Label{1.2}
\no{\eta u}_s\simleq \no{\eta'\dib_b u}_s+\no{\eta'\dib^*_b u}_s+\no{u}_0\quad \T{for any $u\in \mathcal H^\perp\cap C^\infty(b\Om)$}
\\
\T{  in any degree $k\in[0,n-1]$}.
\end{multline}
The same estimate holds for the $\dib$-Neumann problem.
\ep
\bpf
%We start from the so called ``twisted" estimates in the boundary (cf. \cite{S10} Proposition 2.4 p. 18 for the version of these estimates for the $\dib$-Neumann problem).
%We set up tangential basic estimates for a weight adapted to the problem.
We  choose the  orientation $T^\pm$ of the purely imaginary vector field  and consider the  microlocal decomposition of the identity $\T{Id}=\Psi^++\Psi^-+\Psi^0$  and the corresponding decomposition of a form $u=u^++u^-+u^0$ (cf. \cite{K02} Section 2). We recall that $u^0$ enjoys elliptic estimates; we also observe that $[\dib_b^{(*)},\Psi^\pm]\prec\Psi^0$ and hence it suffices to prove \eqref{1.2} separately for $u^+$ and $u^-$. We recall that the star-Hodge operator $u^-\mapsto *\overline{u^-}=*\bar u^+$ 
settles up a correspondence between ``negative" forms in degree $k$ and ``positive" forms in complementary degree $n-1-k$. Thus it suffices to prove \eqref{1.2} for $u^+$. 

We  start from $k\ge1$. We recall the weighted tangential estimates with weight $e^{-\phi}$ for $\phi=-\log\eta^2+t|z'|^2,\,z'\in T^\C_0b\Om$; we point out that even though the weight $\phi$ is not smooth, nevertheless $e^{-\phi}\di_b\phi$ and $e^{-\phi}\di_b\dib_b\phi$ are bounded and hence all integrals below are well defined. Here is the estimate
\begin{equation}
\Label{1.3,5}
\int e^{-\phi}\di_b\dib_b\phi(u^+,u^+)dV+\NO{\bar\nabla u^+}_\phi\simleq\NO{\dib_bu^+}_\phi+\NO{(\dib^*_\phi)_bu^+}_\phi+\NO{u^+}_\phi.
\end{equation}
We first remove $e^{-t|z'|^2}$ from norms since it is uniformely bounded from above and below; thus the  norms in \eqref{1.3,5} change into $\NO{\cdot}_{-\log\eta^2}$. We now describe $\di_b\dib_b\phi$ and $(\dib^*_\phi)_bu^+$. For the first
\begin{equation} 
\Label{1.4,5}
\di_b\dib_b\phi(=\di_b\dib_b(-\log\eta^2+t|z'|^2))=-\frac2\eta\di_b\dib_b\eta+2\frac{\di_b\eta\otimes\dib_b\eta}{\eta^2}+t\di_b z'\otimes \dib_b \bar z'.
\end{equation}
For the second, we start from $(\dib^*_\phi)_b=\dib^*_b+\di_b\log\eta^2-t\bar z'dz'$ and get
\begin{equation}
\Label{1.5,5}
\begin{split}
\NO{(\dib^*_\phi)_bu^+}_\phi&\sim \NO{(\dib^*_{-\log\eta^2+t|z'|^2})_bu^+}_{-\log\eta^2}
\\
&=\NO{\eta\dib_b^*u^+}_0+4\NO{\di_b\eta\contrazione u}_0+\NO{\eta t\bar z'dz'\contrazione u}_0
\\
&+4\Re\Big(\eta\dib_b^* u^+,\di_b\eta\contrazione u^++t\bar z'dz'\contrazione  u^+\Big)_0+4\NO{\di_b\eta\contrazione u^++\eta t\bar z'dz'\contrazione u^+}_0.
\end{split}
\end{equation}
Taking $\T{supp}\,\eta$ in a small neighborhood of $z_o=0$, $t\bar z'$ is also small. By \eqref{1.4,5} and \eqref{1.5,5}, equality \eqref{1.3,5} with $u^+$ replaced by $\Lambda^s u^+$ yields
\begin{equation}
\Label{1.6,5}
\begin{split}
t\NO{\eta u^+}_0+\NO{\eta\bar\nabla u^+}_0&\simleq \NO{\eta\dib_bu^+}_0+\NO{\eta\dib^*_bu^+}_0
\\
&+\NO{\di_b\eta\contrazione u^+}_0+\int\di_b\dib_b(\eta^2)(u^+,u^+)dV.
\end{split}
\end{equation}
Note that, an alternative proof of \eqref{1.6,5} can be obtained from the boundary version of \cite{S10} Proposition 2.4 formula (2.24) with ``twisting coefficient"  $\sqrt a=\eta$ and weight $e^{-t|z'|^2}$.
We apply \eqref{1.6,5} for $u^+$ replaced by $\Lambda^su^+$ and
 wish to do  two operations: to commutate $\dib^{(*)}_b$ with $\Lambda^s$ in the right side of \eqref{1.6,5}, and to estimate the two terms in the second line; (here $\dib^{(*)}_b$ denotes either occurence of $\dib_b$ or $\dib^{*}_b$). For this, we notice that, with the notation $c_s:=\underset {z'}\max|(c_s^j)_j|$, we have
\begin{equation}
\Label{1.4}
\begin{cases}
[\dib^{(*)}_b,\Lambda^s]=(c^j_s)_j\Lambda^s,
\\
\begin{split}
||\di\eta&\contrazione\Lambda^s u^+||^2+\int_{b\Om} \di_b\dib_b\eta^2(\Lambda^su^+,\Lambda^su^+)dV
\\
&
\simleq Q^b_{\eta'\Lambda^{s-\epsilon}}(u^+,u^+)+c_s\NO{\eta'\Lambda^{s-\epsilon}u^+}.
\end{split}
\end{cases}
\end{equation}
Here and in what follows, for an operator $\T{Op}$ such as $\eta'\Lambda^{s-\epsilon}$, we write $Q^b_\T{Op}(u^+,u^+)$ for $\NO{\T{Op}\dib_b u}+\NO{\T{Op}\dib_b^*u}$. 
We use \eqref{1.4} inside \eqref{1.6,5} in which $u^+$ is replaced by $\Lambda^su^+$  and  get
\begin{equation}
\Label{1.5}
\NO{\eta\Lambda^s u^+}_0\simleq Q^b_{\eta\Lambda^s}(u^+,u^+)+Q^b_{\eta'\Lambda^{s-\epsilon}}(u^+,u^+)+c_s\NO{\eta'\Lambda^{s-\epsilon}u^+}_0+\frac{c_s}t \NO{\eta\Lambda^{s}u^+}_0.
\end{equation}
We absorb the term in \eqref{1.5} with a factor of $\frac{c_s}t$ by taking $t$ large  and restart \eqref{1.5} for $\eta$ replaced by $\eta'$ and $\Lambda^s$ by $\Lambda^{s-\epsilon}$ and, by induction on $j$ such that $j\epsilon>s$,  get \eqref{1.2} for any form in degree $1\le k\le n-1$. 

We have to show now that \eqref{1.2} also holds for forms in degree $k=0$. In fact, given $u\in\mathcal H^\perp$, we use that $\dib^*_b$ has closed range, and write
$$
u=\dib^*_bv\quad\T{for some $1$-form $v$ such that $\dib_bv=0$ and $||v||_0\simleq \no{u}_0$.}
$$
We now observe that 
\begin{equation}
\Label{maris}
\begin{cases}
\eta\dib^*_b(v)=\dib^*_b(\eta v)-\di_b\eta\contrazione  v,
\\
\begin{split}
(\dib_b^*v)^+&=\dib^*_bv^+-[\dib^*_b,\Psi^+]v
\\
&=:\dib^*_bv^++v^0,
\end{split}
\end{cases}
\end{equation}
for $v^0\sim-\dot\Psi^+v$.
It follows
\begin{equation}
\Label{1.6}
\begin{split}
||&\eta\Lambda^s u^+||^2=\Big(\Lambda^su^+,\eta^2\Lambda^s\dib^*_bv^+\Big)\underset{\T{denoted $\mathcal E$ below}}{\underbrace{-(\Lambda^su^+,\eta^2\Lambda^sv^0)}}
\\
&=\Big(\Lambda^su^+,\dib^*_b(\eta^2\Lambda^sv^+)\Big)+\Big(\eta\Lambda^su^+,\sum_j(2L_j(\eta)+c_s^j\eta)\Lambda^sv^+_j\Big)+\mathcal E
\\
&=\Big(\eta\Lambda^s\dib_bu^+,\eta\Lambda^sv^+\Big)+2\Re\Big(\eta\Lambda^su^+,\eta\sum_jc_s^j\Lambda^sv^+_j\Big)+\Big(\eta\Lambda^su^+,2\sum_jL_j(\eta)\Lambda^sv^+_j\Big)+\mathcal E
\\
&\underset{\T{\ref{1.5}}}\le\NO{\eta\Lambda^s\dib_b u^+}+sc\NO{\eta\Lambda^su^+}+lc\NO{\eta\Lambda^sv^+}+lc\NO{\di_b(\eta)\contrazione\Lambda^s v^+}+\mathcal E,
\end{split}
\end{equation}
where, sc and lc denote a small and large constant respectively.
Here and in the following, the notation $\mathcal E$ is used for  an error subject to an elliptic gain which can therefore be disregarded.
We have now  to estimate $\NO{\eta\Lambda^sv^+}$ and $\NO{\di_b(\eta)\contrazione\Lambda^s v^+}$. For the second:
\begin{equation}
\Label{contr}
\begin{split}
\NO{\di_b(\eta)\contrazione\Lambda^s v^+}&\simleq Q^b_{\eta'\Lambda^{s-\epsilon}}(v^+,v^+)+c_s\NO{\eta'\Lambda^{s-\epsilon}v^+}
\\
&\simleq \NO{\eta'\Lambda^{s-\epsilon}(\dib^*_bv)^+}+\NO{\eta'\Lambda^{s-\epsilon}v^0}+c_s\NO{\eta'\Lambda^{s-\epsilon}v^+}.
\end{split}
\end{equation}
The central term in the last line above is of type $\mathcal E$.
The first term is $\NO{\eta'\Lambda^{s-\epsilon}u^+}$ which can be controlled by induction. Finally, to handle $\NO{\eta'\Lambda^{s-\epsilon}v^+}$, we apply \eqref{1.5} and get an estimate by means of $\NO{\eta'\Lambda^{s-\epsilon}u^+}+\NO{\eta'\Lambda^{s-\epsilon}v^0}+\NO{\eta''\Lambda^{s-2\epsilon}u^+}+\NO{\eta''\Lambda^{s-2\epsilon}v^+}+\frac{c_s}t\NO{\eta'\Lambda^{s-\epsilon}v^+}$. 
In this way we control $\NO{\di_b(\eta)\contrazione\Lambda^s v^+}$ (and similarly we can control $\int_{b\Om} \di_b\dib_b\eta^2(\Lambda^s v^+,\Lambda^sv^+)dV$).
We pass to $\NO{\eta\Lambda^sv^+}$: since $\dib_b v=0$, then 
\begin{equation}
\Label{extra}
\begin{split}
||&\eta\Lambda^s v^+||^2\underset{\T{\eqref{1.6,5}}}\simleq\frac1t\Big(\NO{\eta\Lambda^s\dib_b^*v^+}+\NO{\di_b(\eta)\contrazione\Lambda^s v^+}+\int_{b\Om} \di_b\dib_b\eta^2(\Lambda^s v^+,\Lambda^sv^+)dV\Big)\\
&\le \frac1t\NO{\eta\Lambda^s u^+}+\NO{\eta\Lambda^sv^0}+\NO{\di_b(\eta)\contrazione\Lambda^s v^+}+\int_{b\Om} \di_b\dib_b\eta^2(\Lambda^s v^+,\Lambda^sv^+)dV.
\end{split}
\end{equation}
In the last line of \eqref{extra}, the first term can be absorbed in the left of \eqref{1.6}, the second is subject to an elliptic gain as $\mathcal E$ above, the third has the estimate \eqref{contr} and the last is similar.  Altogether,  $\NO{\di_b(\eta)\contrazione\Lambda^s v^+}$ and  $\NO{\eta\Lambda^sv^+}$, in \eqref{1.6} are controlled.
Thus induction works in \eqref{1.6} and has the effect of reducing the Sobolev index both of $u^+$ and $v^+$. At the last step, that is at $s=0$, we use the closed range estimate $\no{v^+}_0\le\no{v}_0\simleq \no{u}_0$ and get rid of $v$.

 This concludes the proof of estimate \eqref{1.2} of Proposition~\ref{p1.1}. The corresponding estimate for the $\dib$-Neumann problem is obtained by the technique of \cite{BKZ14}.
\br
We give an alternative proof of \eqref{1.2} for $u^-$ which avoids use of the star-Hodge operator. For this, we start, instead of \eqref{1.3,5} from
\begin{multline}
\Label{concave}
-\int e^\phi\di_b\dib_b\phi(u^-,u^-)dV+\sum_j\int e^\phi\phi_{j\bar j}|u^-|^2dV+\NO{\nabla u^-}_{-\phi}
\\
\simleq \NO{\dib_bu^-}_{-\phi}+\NO{(\dib^*_{-\phi})_bu^-}_{-\phi}+\NO{u^-}_{-\phi}.
\end{multline}
Using the analog of \eqref{1.4,5}, \eqref{1.5,5} with $\phi$ replaced by $-\phi$, we end up with
\begin{multline}
t\NO{\eta u^-}_0+\NO{\eta\nabla u^-}_0\simleq Q_\eta(u^-,u^-)+\NO{\dib_b\eta\wedge u^-}_0
\\
+\int\di_b\dib_b\eta^2(u^-,u^-)dV-2\sum_j\int \eta\eta_{j\bar j}|u^-|^2dV.
\end{multline}
We then use the identities
\begin{equation*}
\begin{cases}
\dib_b\eta\wedge u^-=-*\overline{\di_b\eta\contrazione*\bar u^+},
\\
\di_b\dib_b\eta^2(u^-,u^-)-2\eta\sum_j\eta_{j\bar j}|u^-|^2=\di_b\dib_b\eta^2(*\bar u^+,*\bar u^+),
\end{cases}
\end{equation*}
which shows an action of subelliptic multiplier over $*\bar u^+$. The rest of the proof goes through as before.
\er

\epf

\noindent
{\it Proof of Theorem~\ref{t1.1}. }
By the $L^2$-theory of $\dib_b$, there is well defined in $L^2$ the Green operator $G=\Box_b^{-1}$. As an immediate consequence of \eqref{1.2},
$\dib^*N$ and $\dib N$ have exact local $H^s$-regularity at $z_o$ over $\ker \dib$ and $\ker \dib^*$ respectively. More precisely, we have
\begin{equation}
\Label{1.7}
\no{\eta\dib^{(*)}_bGu}_s\simleq\no{\eta'u}_s+\no{u}_0.
\end{equation}
 Let $S$, resp. $S^*$ be the Szeg\"o, resp. anti-Szeg\"o, projection. By Kohn's formula $S=\T{Id}-\dib^*_bG\dib_b$ and $S^*=\T{Id}-\dib_bG\dib^*_b$, we have that the projections $S^{(*)}$ are also regular, though a loss of one derivative may occur on account of the double application of $\dib^{(*)}_b$. In other words we have
\begin{equation}
\Label{1.8}
\no{\eta S^{(*)}u}_s\simleq \no{\eta'u}_{s+1}+\no{u}_0.
\end{equation}
From this, we can get the (non-exact) regularity of $G$ itself on account of 
\begin{equation*}
\begin{split}
\no{\eta G\alpha}_s&=\no{\eta\Box_bG^2\alpha}_s
\\
&\simleq \no{\eta\dib_b\dib^*_bG^2S\alpha}_s+\no{\eta\dib^*_b\dib_bG^2S^*\alpha}
\\
&\simleq \no{\eta\dib_bG\dib^*_bGS\alpha}_s+\no{\eta\dib^*_bG\dib_bGS^*\alpha}_s
\\
&\simleq \no{\eta'\dib^*_bGS\alpha}_s+\no{\eta'\dib_bGS^*\alpha}_s
\\
&\simleq \no{\eta''S\alpha}_s+\no{\eta''S^*\alpha}_s
\\
&\simleq \no{\eta'''\alpha}_{s+1}.
\end{split}
\end{equation*}
This estimate with loss of 1 derivative is an ``a-priori" estimate. The method of the elliptic regularization makes it a ``genuine" estimate; this clearly suffices for local $C^\infty$-regularity of the Green operator $G$.  The similar conclusion on the $C^\infty$-regularity of the Neumann operator $N$ is obtained from the variant of \eqref{1.2} for the $\dib$-Neumann problem.

\hskip12cm $\Box$

\section{A class of Examples}
A large class of domains to which Theorem~\ref{t1.1} applies is provided by the following
\bt
\Label{p2.1}
In $\C^n$ we consider a ``decoupled" pseudoconvex domain whose boundary $b\Om$ is defined in a neighborhood of $0$ by
$$
2x_n=\sum_{j=1}^{n-1}h^j(z_j),
$$
for $h^j$ real subharmonic, that is, satisfying $h^j_{j\bar j}\ge0$. We make the additional assumptions that each $h^j$ has finite type $2m_j$ for $z_j\neq0$ and that, up a harmonic term $\Re F^j$, we have $|h^j_j+\Re F^j|\simleq h_{j\bar j}^j$.

Then, for a fundamental system of cut-off $\eta$ at $0$,  $\di_b\eta$ and $\di_b\dib_b\eta^2$ are $\frac1{2m}$-subelliptic multipliers for $m=\sup_jm_j\ge2$ over forms $u^+$ in degree $k\in [0,n-1]$. 
\et
\bpf
We choose  a cut-off $\chi$ in $\R$ at $0$, set $\zeta=\Pi_j\chi(|z_j|)$, $\theta=\chi(y_n)$, and define $\eta=\zeta\theta$.
We also write a general coefficient of $u$ in degree $k$ as $u_{jK}$ for $j=1,...,n-1$ and $|K|=k-1$; we also use the 
 notation $r:=2x_n-\sum_jh^j$. The crucial point in the proof below is that, $r$ being decoupled, we have
\begin{equation}
\Label{3.2}
\sumK\sumij r_{i\bar j}u_{iK}\bar u_{jK}-\sumj\sumK r_{j\bar j}|u_{jK}|^2=0.
\end{equation}
Thus,  the basic estimate not only yields
\begin{equation}
\Label{3.1}
\sum_j\sumJ\NO{\bar L_ju^+_J}_0\simleq Q^b(u^+,u^+)+\NO{u^+}_0,
\end{equation}
as usual, but also
\begin{equation}
\Label{3.3}
\sum_j\sumK\NO{ L_ju_{jK}^+}_0\simleq Q^b(u^+,u^+)+\NO{u^+}_0.
\end{equation}
We select an index $j_o$. 
Since, the iterated brackets $\underset{2m_{j_o}}{\underbrace{[\overset{(-)}L_{j_o},[\overset{(-)}L_{j_o},[...]]]}}$ (where $\overset{(-)}L_{j_o}$ denotes either occurence of $L_{j_o}$ or $\bar L_{j_o}$) generate the purely imaginary vector field $T=\di_{y_n}$ over $\T{supp}\,\dot\zeta_{z_{j_o}}\subset\{z_{j_o}:\,z_{j_o}\neq0\}$, then we have
\begin{equation}
\Label{3.4}
\begin{split}
\no{T(\zeta_{z_{j_o}}u^+_{j_oK})}_{-1+\frac1{2m_{j_o}}}&\simleq\no{L_{j_o}(\zeta_{z_{j_o}}u^+_{j_oK})}_0+\no{\bar L_{j_o}(\zeta_{z_{j_o}}u^+_{j_oK})}_0+\no{u^+_{j_oK}}_0\\
&\underset{\T{\eqref{3.1}, \eqref{3.3}}}\simleq Q^b(u^+,u^+)+\NO{u^+_{j_oK}}_0.
\end{split}
\end{equation}
Thus $Q^b+\NO{\cdot}_0$ contains, over $\T{supp}\,\zeta_{z_{j_o}}$, the norm of a fractional derivative $\NO{T^{\frac1{2m_{j_o}}}u^+_{j_oK}}_0$  and of a full derivative $\NO{\bar L_ju^+_{j_oK}}_0$ for any $j=1,...,n-1$. As for  $L_j$, this is  already contained in $Q^b+\NO{\cdot}_0$ for $j=j_o$ according to \eqref{3.3}. For $j\neq j_o$, we have to change  $L_j$  into $\bar L_j$. For this, we use the identity
$$
\NO{L_j(\zeta_{z_{j_o}} u_{j_oK}^+)}_0\simleq\NO{\bar L_j(\zeta_{z_{j_o}} u_{j_oK}^+)}_0+\Big([L_j,\bar L_j](\zeta_{z_{j_o}}u^+_{j_oK}),\zeta_{z_{j_o}}u^+_{j_oK}\Big)+\NO{u^+_{j_oK}}_0;
$$
next, we express the commutator as $[L_j,\bar L_j]=r_{j\bar j}T+\sum_h\overset{(-)}a_h\overset{(-)}L_h$. The terms $\Big( \zeta_{z_{j_o}}\overset{(-)}a_h\overset{(-)}L_hu^+_{j_oK}, \zeta_{z_{j_o}}u^+_{j_oK}\Big)$  can be estimated by $sc\sum_j\NO{ \zeta_{z_{j_o}}\overset{(-)}L_ju^+_{j_oK}}_0+lc\NO{u^+_{j_oK}}_0$ which yields
\begin{equation}
\Label{3.5}
\begin{split}
\sum_j&\no{L_j(\zeta_{z_{j_o}}u^+_{j_oK})}_{-1+\frac1{2m_{j_o}}}^2\simleq \sum_j\no{\bar L_j(\zeta_{z_{j_o}}u^+_{j_oK})}_{-1+\frac1{2m_{j_o}}}^2+
 \no{T(\zeta_{z_{j_o}}u^+_{j_oK})}_{-1+\frac1{2m_{j_o}}}^2
\\
&+sc\Big(\sum_j\no{ L_j(\zeta_{z_{j_o}}u^+_{j_oK})}^2_{-1+\frac1{2m_{j_o}}}
+\sum_j\no{ \bar L_j(\zeta_{z_{j_o}}u^+_{j_oK})}_{-1+\frac1{2m_{j_o}}}^2\Big)+lc\NO{u^+_{j_oK}}_0
\\
&\underset{\T{\eqref{3.4}}}\simleq Q^b(u^+,u^+)+\underset{\T{absorbed}}{\underbrace{
sc\sum_j\no{ L_j(\zeta_{z_{j_o}}u^+_{j_oK})}_{-1+\frac1{2m_{j_o}}}^2}}+lc\NO{u^+_{j_oK}}_0.
\end{split}
\end{equation}
Taking summation over $j_o$ and $K$, and the minimum $\frac1{2m}$ of the $\frac1{2m_{j_o}}$'s, we get the estimate for $\di_b\zeta=\di_b\Pi_j\zeta_j$
\begin{equation}
\Label{2.3}
\begin{split}
\NO{\di_b\zeta\contrazione u^+}_{\frac1{2m}}&\simleq Q^b(u^+,u^+)+\NO{u^+}_0
\\
&\simleq Q^b(u^+,u^+),
\end{split}
\end{equation}
where the second estimate follows from the closed range.
Passing to a general $\di_b\eta=\di_b\zeta\theta$,
we notice that
\begin{equation}
\Label{2.4}
\di_b\eta=(L_j\eta)_{j=1,...,n-1}=(\underset{\T{(a)}}{\underbrace{\zeta_{z_j}\theta}})_{j=1,...,n-1}+(\underset{\T{(b)}}{\underbrace{\zeta h^j_{z_j}\dot\theta}})_{j=1,...,n-1}.
\end{equation}
Now, (a) has already been estimated in \eqref{2.3}. As for (b), we observe that in new complex coordinates in which we get rid of harmonic terms in the $h^j$'s, we have by hypothesis $|h^j_{z_j}|^2\simleq h^j_{z_j\,\bar z_j}$. It follows 
\begin{equation}
\begin{split}
\NO{(h^j_{z_j})_{j=1,...,n-1}\contrazione u^+}_{\frac12}&\simleq \int \di_b\dib_br(T^{\frac12}u^+,T^{\frac12}u^+)dV+\NO{u^+}_0
\\
&\simleq Q^b(u^+,u^+);
\end{split}
\end{equation}
again, we have estimated $\NO{u^+}_0\simleq Q^b$ by closed range.
This, together with \eqref{2.3}, shows that the gradient $\di_b\eta$ is a $\frac1{2m}$-subelliptic multiplier. We pass to the Levi form. We start from the obvious equality $L_i\bar L_j(\eta^2)=2L_i(\eta)\bar L_j(\eta)+2\eta L_i\bar L_j(\eta)$ and
$$
L_i\bar L_j(\eta)=(\zeta_{z_i}\zeta_{\bar z_j}\theta+\zeta_{z_i\bar z_j}\theta)+(\zeta_{z_i}h^j_{\bar z_j}\dot\theta+\zeta_{\bar z_j}h^i_{z_i}\dot\theta)+\zeta h^j_{z_i\,\bar z_j}\dot\theta+\zeta h^i_{z_i}h^j_{\bar z_j}\ddot\theta.
$$
Now, the first and second  terms in the right are controlled by (a) of \eqref{2.4} above. The  third and fourth  by 
$$
\Big|\int \zeta_{z_i}u^+_ih^j_{\bar z_j}\bar u^+_j\dot\theta dV\Big|\underset{\T{Cauchy-Schwarz}}\le \no{\zeta_{z_i}u^+_i}_0\,\no{h^j_{\bar z_j}\bar u^+_j\dot\theta}_0,
$$
and then by (a) combined with (b). The fifth by
\begin{equation*}
\begin{split}
\sum_j\NO{ T^{\frac12}(h^j_{z_j\bar z_j}u^+_j)}_0&\simleq \int \di_b\dib_br( T^{\frac12}u^+, T^{\frac12}u^+)dV+\NO{u^+}_0\\
&\simleq Q^b(u^+,u^+).
\end{split}
\end{equation*}
Finally, the sixth by (b).

\epf

\be
For the pseudoconvex domain with boundary defined, in a neighborhood of $0$, by
$$
2x_n=\sum_{j=1}^{n-1}e^{-\frac1{|z_j|^{s_j}}}x_j^{2m_j}\quad\T{any $s_j>0$ and $m_j\ge0$},
$$
we can readily verify that the hypotheses of Proposition~\ref{p2.1} are satisfied. Hence, on account of Theorem~\ref{t1.1}, $\Box_b$ and $\Box$ are $C^\infty$-hypoelliptic at $0$.
\ee


\begin{thebibliography}{BKZ09}
%\bibitem{B13}{\bf L. Baracco}---The range of the tangential Cauchy Riemann system to a $CR$- embedded manifold, {\em  Invent. Math. } {\bf 190} (2012), 505--510
%\bibitem{BKP}{\bf L. Baracco, T.V. Khanh and S. Pinton}---Uniform regularity in a wedge and regularity of traces of $CR$ functions, {\em  J. Geom. Anal.} {\bf 20} (2010), 996--1007
\bibitem{BKZ12} {\bf L. Baracco, T.V. Khanh and G. Zampieri}---Propagation of regularity for solutions of the Kohn Laplacian in a flat boundary, {\em Adv. Math.} {\bf 230} (2012) 1972--1978
\bibitem{BKZ14} {\bf L. Baracco, T.V. Khanh and G. Zampieri}---Hypoellipticity of the $\bar\partial$-Neumann problem at a point of infinite type, {\em Asian J. Math.} (2014)
\bibitem{BPZ13} {\bf L. Baracco, S. Pinton and G. Zampieri}---Hypoellipticity of the $\dib$-Neumann problem at a set of infinite type with positive CR dimension, (2013)
\bibitem{C87} {\bf D. Catlin}---Subelliptic estimates for the $\bar{\partial}$-Neumann problem on pseudoconvex domains, {\em  Ann. of Math. } {\bf 126} (1987), 131-191
\bibitem{Ch00}{\bf M. Christ}---Hypoellipticity: geometrization and speculation, {\em Progress in Math. Birkh"auser Basel},    {\bf188} (2000), 91--109
\bibitem{Ch02} {\bf M. Christ}---Hypoellipticity of the Kohn Laplacian for three-dimensional tubular Cauchy-Riemann structures, {\em J. of the Inst. of Math. Jussieu} {\bf 1} (2002), 279--291
\bibitem{DT95}{\bf M. Derridj and D.S. Tartakoff}---Microlocal analiticity for $\Box_b$ in block-decoupled pseudoconvex domains, {\em Math. Z.} {\bf 220} (1995), 477--493
%\end{document}
%\bibitem{CS01} {\bf S.C. Chen and  M.C. Shaw}--Partial differential equations in several complex variables, {\em Studies in Adv. Math. - AMS Int. Press} {\bf 19} (2001)
%\bibitem{F71}{\bf V.S. Fedii}---On a criterion for hypoellipticity, {\em Math. USSR Sb.}, {\bf 14} (1971), 15--45
%\bibitem{E97} {\bf L. Evans}---Partial Differential Equations, {\em Graduate Studies in Math.} {\bf 19} (1997)
\bibitem{FK72}{\bf G.B. Folland and J.J. Kohn}---The Neumann problem for the Cauchy-Riemann complex, {\em Ann. Math. Studies, Princeton Univ. Press, Princeton N.J.} {\bf 75} (1972)
\bibitem{HT83}{\bf J. Hanges and F. Treves}---Propagation of holomorphic extandibility of CR functions, {\em Math. Ann.} {\bf 263} n. 2 (1983), 157--177
%\bibitem{D82}{\bf J. D'Angelo}---Real hypersurfaces, order of contact, and applications, {\em Ann. of Math.} {\bf 115} (1982), 615--637
%\bibitem{K64} {\bf J.J. Kohn}---Harmonic integrals on strongly pseudoconvex manifolds, I, {\em Ann. Math.} {\bf 78} (1963), 112--148; II, {\em Ann. Math.} {\bf 79} (1964), 450--472 
%\bibitem{K72}{\bf J.J. Kohn}---Boundary behavior of $\bar\partial$ on weakly pseudo-convex manifolds of dimension two, {\em J. Diff. Geom.} {\bf 6} (1972), 523--542
%\bibitem{K73} {\bf J.J. Kohn}---Global regularity for $\bar\partial$ on weakly pseudo-convex manifolds, {\em Trans. of the A.M.S.} {\bf 181} (1973), 273--292
%\bibitem{K77} {\bf J.J. Kohn}---Methods of partial differential equations in complex analysis,{\em Proceedings of Symposia in pure Mathematics} {\bf 30}(1977), 215--237
%\bibitem{Ke72} {\bf N. Kerzman}---The Bergman kernel function. Differentiability at the boundary, {\em Math. Ann.} {\bf 195} (1972), 149--158
%\bibitem{Kh09} {\bf T.V. Khanh}---A general method of weights in the $\bar\partial$-Neumann problem, Ph.D. Thesis, Padua (2009)
%\bibitem{KZtan09}{\bf T.V. Khanh and G. Zampieri}---Estimates for regularity of the tangential $\bar\partial$-system, {\em Math. Nachr.} {\bf 284}  no. 17-18 (2011), 2212-2224
%\bibitem{KPZ12}{\bf T.V. Khanh, S. Pinton and G. Zampieri}---Compactness estimates for $\Box_b$ on a CR
%manifold, {\em Proc. of the A.M.S.} {\bf 140} (2012),  3229 − 3236
\bibitem{KZ10}{\bf T.V. Khanh and G. Zampieri}---Regularity of the $\bar\partial$-Neumann problem at a flat point, {\em J. Funct. Anal.} {\bf 259} no. 11 (2010), 2760-2775
%\bibitem12}{\bf T.V. Khanh and G. Zampieri}---Necessary geometric and analytic conditions for a general estimate in the $\dib$-Neumann problem, {\em Invent. Math.}  {\bf 188} (2012), 729--750.
\bibitem{K79} {\bf J.J. Kohn}---Subellipticity of the $\bar\partial$-Neumann problem on pseudoconvex domains: sufficient conditions, {\em Acta Math.} {\bf 142} (1979), 79--122
\bibitem{K00} {\bf J.J. Kohn}---Hypoellipticity at points of infinite type, {\em Contemporary Math.} {\bf 251} (2000), 393--398
\bibitem{K02} {\bf J.J. Kohn}---Superlogarithmic estimates on pseudoconvex domains and CR manifolds, {\em Annals of Math.} 
{\bf 156} (2002), 213--248
\bibitem{KN65} {\bf J.J. Kohn and L. Nirenberg}---Non-coercive boundary value problems, {\em Comm. Pure Appl. Math.} {\bf 18} (1965), 443--492
%\bibitem{MN92} {\bf  J. D. McNeal}--- Lower bounds on the Bergman metric near a point of finite type.  {\em Ann. of Math. }   136  (1992),  {\bf  2,} 339--360.
\bibitem{S10} {\bf E. Straube}---Lectures on the $L^2$-Sobolev theory of the $\bar\partial$-Neumann problem, {\em ESI Lect. in Math. and Physics} (2010)
%\bibitem{Z08} {\bf G. Zampieri}---Complex Analysis and CR Geometry, {\em AMS ULECT} {\bf 43} (2008)
%\bibitem{H65} {\bf L. Hormander}---$L^2$ estimates and existence theorems for the $\bar{\partial}$ operator, {\em Acta Math.} {\bf 113} (1965), 89--152
%\bibitem{H73}{\bf L. H\"ormander}---An intoduction to complex analysis in several complex variables, {\em  Van Nostrand} Princeton N.J. (1973)
%\bibitem{Ho85} {\bf L.H. Ho}---Subellipticity of the $\bar{\partial}$-Neumman problem on the nonpseudoconvex, {\em Trans. AMS} {\bf 291} (1985), 43-73 

%\bibitem{Ho91} {\bf L.H. Ho}---Subellipticity of the $\bar{\partial}$-Neumman problem for $n-1$ forms, {\em Trans. AMS} {\bf 325} (1991), 171-185 
%\bibitem{M87} {\bf Y. Morimoto}---A criterion for hypoellipticity of second order differential operators, {\em Osaka J. Math.}, {\bf 24} (1987), 651--675

%\bibitem{TZ1}{\bf T.V. Khanh and G. Zampieri}---Subellipticity of the $\bar\partial$-Neumann problem on a weakly q-pseudoconvex/concave domain, {\em arXiv:0804.3112v} (2008)
%\bibitem{TZ2}{\bf T.V. Khanh and G. Zampieri}---Compactness of the $\bar\partial$-Neumann operator on a q-pseudoconvex domain, (2008)
%\bibitem{Z08} {\bf G. Zampieri}---Complex analysis and CR geometry, {\em AMS ULECT} {\bf 43} (2008)
\end{thebibliography}
 \end{document}